\documentclass[11pt]{amsart}
\usepackage{amssymb,latexsym,amsmath,amsthm,amsfonts,amsxtra}

\theoremstyle{plain}
\newtheorem{theorem}{Theorem}
\newtheorem{corollary}{Corollary}

\newtheorem{lemma}{Lemma}

\theoremstyle{definition}
\newtheorem{definition}{Definition}

\theoremstyle{remark}

\numberwithin{equation}{section}

\def\R{{\mathbb{R}}}
\def\CC{{\mathbb{C}}}

\def\Z{{\mathbb{Z}}}

\def\dir{{\hbox{dir}}}

\begin{document}

\title{Uniform formulae for coefficients of meromorphic functions in two variables. Part I.}

\author{Manuel Lladser}

\address{Department of Applied Mathematics\\
University of Colorado\\
PO Box 526 UCB\\
Boulder, CO 80309-0526\\
USA.}

\email{manuel.lladser@colorado.edu}

\urladdr{http://amath.colorado.edu/faculty/lladser}



\date{Version of April 05 2006.}

\begin{abstract}
Uniform asymptotic formulae for arrays of complex numbers of the form $(f_{r,s})$, with $r$ and $s$ nonnegative integers, are provided as $r$ and $s$ converge to infinity at a comparable rate.
Our analysis is restricted to the case in which the generating function $F(z,w):=\sum f_{r,s} z^r w^s$ is meromorphic in a neighborhood of the origin. We provide uniform asymptotic formulae for the coefficients $f_{r,s}$ along directions in the $(r,s)$-lattice determined by regular points of the singular variety of $F$. Our main result derives from the analysis of a one dimensional parameter-varying integral describing the asymptotic behavior of $f_{r,s}$. We specifically consider the case in which the phase term of this integral has a unique stationary point, however, allowing the possibility that one or more stationary points of the amplitude term coalesce with this. Our results find direct application in certain problems associated to the Lagrange inversion formula as well as bivariate generating functions of the form $v(z)/(1-w\cdot u(z))$.
\end{abstract}

\maketitle

\section{Introduction}

Suppose that $G(z,w)$ and $H(z,w)$ are analytic functions of the complex variables $z$ and $w$ in an open polydisk centered at the origin and assume that $H(0,0)\ne0$. Then, the function
\[F(z,w):=\frac{G(z,w)}{H(z,w)}\]
is analytic in a neighborhood of the origin in $\CC^2$; in particular, it has a power series expansion of the form $\sum f_{r,s}z^rw^s$, where the indices $r$ and $s$ are nonnegative integers. In what follows we use the notation $[z^rw^s]F$ to refer to the coefficient of $z^rw^s$ in the power series expansion of $F$. We also use the notation $(r,s)\to\infty$ as a shorthand for $r\to\infty$ and $s\to\infty$.

Generating functions of the above form occur frequently in the study of discrete random structures and analysis of algorithms (see~\cite{PW05} for a comprehensive account of examples). For a wide class of bivariate functions of this kind the coefficients $[z^rw^s]F$ are expected, up to an exponential factor, to be of order $s^{-(p+1)/n}$ as $(r,s)\to\infty$ with $r/s$ fixed. Here, the coefficients $p$ and $n$ are functions of the ratio $r/s$. In particular, the asymptotic behavior of $[z^rw^s]F$ can be understood even if $r/s$ varies but in such a way that $p$ and $n$ do not change. In this paper we show how to provide uniform asymptotic formulae for the coefficients $[z^rw^s]F$ as $(r,s)\to\infty$ when $r/s$ is restricted to a set of values where the coefficient $p$ may not remain constant.

From this point on we assume as given a point $(\zeta,\omega)$ which is a strictly minimal simple zero of $H$. By simple zero we mean that $H(\zeta,\omega)=0$, however, the complex gradient $\nabla H(\zeta,\omega)\ne0$. By strictly minimal zero we mean that $\zeta\cdot\omega\ne0$ and that $(\zeta,\omega)$ is the only zero of $H(z,w)$ in the polydisk where $|z|\le|\zeta|$ and $|w|\le|\omega|$.

The pioneering work of Pemantle and Wilson~\cite{PW01} implies that
it is possible to determine an asymptotic expansion for the coefficients of $F$ only along certain direction in the $(r,s)$-lattice specified by $(\zeta,\omega)$. This direction corresponds to the line
\begin{equation}\label{def:dir(z0,w0)}
\dir(\zeta,\omega):=\{(r,s)\in\R^2: r\cdot\omega\,H_w(\zeta,\omega)-s\cdot\zeta\,H_z(\zeta,\omega)=0\}\,,
\end{equation}
where $H_w$ and $H_z$ respectively denote the complex partial derivative of $H$ with respect to $w$ and $z$. For simplicity it will be assumed ahead that $H_w(\zeta,\omega)\ne0$. In particular, $(r,s)\in\dir(\zeta,\omega)$ if and only if $r/s=d(\zeta,\omega)$, where
\[d(\zeta,\omega):=\frac{\zeta\,H_z(\zeta,\omega)}{\omega\,H_w(\zeta,\omega)}\,.\]
The strict minimality of $(\zeta,\omega)$ implies that $d(\zeta,\omega)\ge0$ (see Lemma 2.1 in~\cite{PW01}). Furthermore, if this quantity is a rational number, Pemantle and Wilson show that there are integers $n=n(\zeta,\omega)\ge2$ and $p=p(\zeta,\omega)\ge0$ and coefficients $c_j=c_j(\zeta,\omega)$, with $j\ge p$ and $c_p\ne0$, such that
\begin{equation}\label{int:[zrws] using (z0,w0)}
[z^rw^s]F\approx\frac{\zeta^{-r}\omega^{-s}}{2\pi}\sum_{j=p}^\infty c_j\,s^{-(j+1)/n}\,,
\end{equation}
for all $(r,s)\in\dir(\zeta,\omega)$, as $(r,s)\to\infty$. The asymptotic notation used above is in the standard sense where the sequence $(s^{-(j+1)/n})_{j\ge p}$ is the so called auxiliary asymptotic sequence. This means that the difference between the left and the right-hand side term above, with the summation truncated to the term in which $j=m$, is $O(\zeta^{-r}\,\omega^{-s}\,s^{-(m+2)/n})$, as $(r,s)\to\infty$. 

The technique used to obtain the asymptotic formula in (\ref{int:[zrws] using (z0,w0)}) proceeds by relating the coefficient $[z^rw^s]F$, with $(r,s)\in\hbox{dir}(\zeta,\omega)$, to a one dimensional Fourier or Laplace like integral of the form
\begin{equation}\label{int:FL-integral}
\frac{\zeta^{-r}\omega^{-s}}{2\pi}\int a(\theta;\zeta,\omega)\exp\{-s\cdot f(\theta;\zeta,\omega)\}d\theta\,.
\end{equation}
We will refer loosely to $a(\theta;\zeta,\omega)$ and $f(\theta;\zeta,\omega)$ respectively as the derived amplitude and phase term. Roughly speaking, the expansion in (\ref{int:[zrws] using (z0,w0)}) is in powers of $s^{-1/n}$ because the derived phase term vanishes to degree $n$ in the variable $\theta$ about $\theta=0$, which turns out to be the dominant critical point of the integral. Furthermore, $[z^rw^s]F$ is of order $\zeta^{-r} \,\omega^{-s}\,s^{-(p+1)/n}$ because the derived amplitude term vanishes to degree $p$ at $\theta=0$. (See~\cite{Bru81} and~\cite{BleHan86} for a compelling introduction to the main techniques used to study the asymptotic behavior of Fourier-Laplace integrals.)

The asymptotic expansion in (\ref{int:[zrws] using (z0,w0)}) holds usually along a wider set of directions in the $(r,s)$-lattice. Indeed, suppose that $K$ is a compact set of strictly minimal simple zeros of $H$ and consider the set
\begin{equation}\label{int:Lambda}
\Lambda:=\bigcup\limits_{(\zeta,\omega)\in K}\left\{(r,s)\in\R^2:\,\frac{r}{s}=d(\zeta,\omega)\right\}\,.
\end{equation}
Observe that $\Lambda$ is a cone if $K$ is connected.
Theorem 3.3 in~\cite{PW01} implies that (\ref{int:[zrws] using (z0,w0)}) holds uniformly as $(r,s)\to\infty$, with $(r,s)\in\Lambda$, provided that the derived amplitude and phase term in~(\ref{int:FL-integral}) do not change their degree of vanishing about $\theta=0$ as $(\zeta,\omega)$ varies over $K$. In particular, if for each $(r,s)\in\Lambda$, $(\zeta,\omega)=(\zeta(r,s),\omega(r,s))\in K$ is such that $d(\zeta,\omega)=r/s$ then
\begin{equation}\label{int:[zrws]F leading order}
[z^rw^s]F\sim c_p(\zeta,\omega)\cdot\frac{\zeta^{-r}\omega^{-s}}{2\pi}\cdot s^{-(p+1)/n}\,,
\end{equation}
uniformly for all $(r,s)\in\Lambda$, as $(r,s)\to\infty$. The notation used above means that the ratio of the two sides tends to $1$ as $(r,s)\to\infty$. 

However, when there is a change of degree of the derived amplitude or phase term in (\ref{int:FL-integral}) the hypotheses of~\cite{PW01} are not met and therefore no conclusion may be drawn from it. When the change of degree is in the phase term one needs to build a bridge between differently scaled regions. This is hard work and will be presented in the forthcoming paper~\cite{Lla05}. (For further details on this case see Theorem 6.6 in Chapter 6 in~\cite{Lla03}.)

The main contribution of the present paper is to settle the case in which only the derived amplitude term may undergo a change of degree. Although it is not mentioned in~\cite{PW01} the asymptotic formula in (\ref{int:[zrws] using (z0,w0)}) is still valid for $(r,s)\in\Lambda$ but it requires a more careful interpretation. To amplify on this consider the case in which at a particular point $(\zeta_c,\omega_c)\in K$, the derived amplitude term in (\ref{int:FL-integral}) vanishes to degree $q$ yet, for all $(\zeta,\omega)\in K$ nearby $(\zeta_c,\omega_c)$, the derived amplitudes vanish to some degree $p<q$. Then (\ref{int:[zrws]F leading order}) implies that up to the exponential factor $\zeta^{-r}\omega^{-s}$,
\begin{equation}\label{int:leading order of [zrws]F}
[z^rw^s]F\hbox{ is of order }\left\{\begin{array}{l}
s^{-(q+1)/n},\text{ if }r/s=d(\zeta_c,\omega_c);\\
s^{-(p+1)/n},\text{ otherwise },
\end{array}\right.
\end{equation}
as $(r,s)\in\Lambda$ goes to infinity, provided that $r/s$ remains constant.

A problem of interest is how to bridge the gap between the asymptotic orders in (\ref{int:leading order of [zrws]F}) as $r/s$ approaches $d(\zeta_c,\omega_c)$, as $(r,s)\to\infty$. As we shall see in the coming section, we resolve this problem with great generality, and we can provide a uniform asymptotic expansion for the coefficients $[z^rw^s]F$ as long as $(r,s)\in\Lambda$ and $|r/s-d(\zeta_c,\omega_c)|$ is sufficiently small. Our main result builds on the asymptotic analysis of an integral such as the one in~(\ref{int:FL-integral}) which does not rely on having the term $a(\theta;\zeta,\omega)$ vanish to constant degree as $(\zeta,\omega)$ varies over $K$. The technique we propose to analyze integrals of this kind draws on the techniques of Chester et al.~\cite{CFU57}, the results of Levinson on polynomial canonical representations~\cite{Lev60I},~\cite{Lev60II}, and the work of Pemantle and Wilson~\cite{PW01}. All these techniques are founded on complex variable methods. For a compelling introduction to function theory of one 
or several complex variables see~\cite{Cartan73},~\cite{Rudin87} or~\cite{Taylor02}.

Under appropriate hypotheses, our main result implies that $[z^rw^s]F$ has (up to an exponentially decreasing factor) an asymptotic expansion of the form
\begin{equation}\label{int:series [zrws]F}
[z^rw^s]F=\sum_{j=p}^q c_j(r/s)\cdot s^{-(j+1)/n}+o(s^{-(q+1)/n})\,,
\end{equation}
where the coefficients $c_j(r/s)$ are analytic functions of $r/s$ and, except for $j=q$, they all vanish when $r/s=d(\zeta_c,\omega_c)$.  Furthermore, the above expansion is uniform as $(r,s)\to\infty$ provided that $r/s$ is sufficiently close to $d(\zeta_c,\omega_c)$. Observe how the condition of having $c_j(d(\zeta_c,\omega_c))=0$, for $j\ne q$, and $c_q(d(\zeta_c,\omega_c))\ne0$ explains the asymptotic behavior described in~(\ref{int:leading order of [zrws]F}).

In well behaved situations one finds that the coefficients $c_j(r/s)$ in (\ref{int:series [zrws]F}) are all nonnegative. This sign constraint prevents cancellation between the terms participating in the summation in~(\ref{int:series [zrws]F}). As a result, one obtains that
\[[z^rw^s]F=(1+o(1))\cdot\sum_{j=p}^q c_j(r/s)\cdot s^{-(j+1)/n}\,,\]
if $r/s\to d(\zeta_c,\omega_c)$ as $(r,s)\to\infty$. On the contrary, when the coefficients $c_j(r/s)$ have mixed signs, and depending on the rate at which $r/s$ approaches to $d(\zeta_c,\omega_c)$, it is possible that terms in the summation in~(\ref{int:series
[zrws]F}) cancel one another and therefore
\[\sum_{j=p}^q c_j(r/s)\cdot s^{-(j+1)/n}=o(s^{-(q+1)/n})\,.\]
Thus, in situations where $r/s\to d(\zeta_c,\omega_c)$ in such away that the above asymptotic formula holds, our main result allow us to conclude only that $[z^rw^s]F$ is $o(s^{-(q+1)/n})$ as $(r,s)\to\infty$. The effect of cancellation to determine asymptotic formulae for the coefficients $[z^rw^s]F$ is illustrated in Application 2 on the next section.

\section{Main Definitions and Results with Applications}

To state our main definition we recall that if $U(z,w)$ is analytic in an open neighborhood of a point $(z_0,w_0)$ in $\CC\times\CC$ then it is possible to represent $U$ in the form
\[U(z,w)=\sum_{k=0}^\infty U_k(z)\cdot(w-w_0)^k\,,\]
where $U_k(z):=\frac{\partial^k U}{\partial w^k}(z,w_0)$. The above series is usually referred to as the Hartog's series of $U$ in powers of $(w-w_0)$ about the point $(z_0,w_0)$. This series is uniformly convergent for all $(z,w)$ 
in polydisks of the form $\{(z,w): |z-z_0|\le\epsilon, |w-w_0|\le\epsilon|\}$ provided that the polydisk is completely contained in the domain where $U$ is analytic (see Section 4.5 in~\cite{Lla03}).

To state our main result the following definition will be used.

\begin{definition}
Given nonnegative integers $p<q$ and a function $U(z,w)$ analytic in an open neighborhood of a point $(z_0,w_0)$ in $\CC\times\CC$, we say that $U$ has a $p$-to-$q$ change of degree about $w=w_0$ as $z\to z_0$ provided that the Hartog's series of $U$ in powers of $(w-w_0)$ about the point $(z_0,w_0)$ is of the form $U(z,w)=U_p(z)\cdot(w-w_0)^p+\ldots+U_q(z)\cdot(w-w_0)^q+\ldots$ where $U_j(z_0)=0$, for all $p\le j<q$, however, $U_q(z_0)\ne0$.  On the contrary, if $U(z,w)=U_p(z)\cdot(w-w_0)^p+\ldots$ with $U_p(z_0)\ne0$ we say that $U$ vanishes to constant degree $p$ about $w=w_0$ as $z\to z_0$. Alternatively, we will sometimes say that $U$ has a $p$-to-$p$ change of degree about $w=w_0$ as $z\to z_0$.
\end{definition}

In what follows, $G(z,w)$ and $H(z,w)$ are given analytic functions in some open polydisk $D$ centered at the origin in $\CC^2$ and it is assumed that $H(0,0)\ne0$. We also assume as given a compact set $K\subset D$ of strictly minimal simple zeros of $H$ containing a particular point $(\zeta_c,\omega_c)$ such that $H_w(\zeta_c,\omega_c)\ne0$. The Implicit Function Theorem (see IV.5.6 in~\cite{Cartan73}) lets us then parametrize the zero set of $H$ near $(\zeta_c,\omega_c)$ in the form $\omega=g(\zeta)$, where $g$ is certain analytic function of $\zeta$ near $\zeta=\zeta_c$.

For each $(\zeta,\omega)\in K$, $\dir(\zeta,\omega)$ is the line defined as in~(\ref{def:dir(z0,w0)}) and $\Lambda$ is the cone defined in (\ref{int:Lambda}). For each $(r,s)\in\Lambda$ such that $r/s=d(\zeta_c,\omega_c)$ we define $(\zeta(r,s),\omega(r,s)):=(\zeta_c,\omega_c)$. Furthermore, for each $(r,s)\in\Lambda$ we let $(\zeta,\omega)=(\zeta(r,s),\omega(r,s))\in K$ be such that $(r,s)\in\dir(\zeta,\omega)$. For the validity of our main result, we require the continuity condition
\[(\zeta(r,s),\omega(r,s))\to(\zeta_c,\omega_c)\,,\]
as $r/s\to d(\zeta_c,\omega_c)$. Indeed, since
\[ (r,s)\in\dir(\zeta,\omega)\iff\frac{r}{s}=-\frac{\zeta g'(\zeta)}{g(\zeta)}\,,\]
to satisfy the continuity condition it is enough to select $\zeta(r,s)=\zeta$ and $\omega(r,s)=g(\zeta)$, where $\zeta$ is the closest solution to $\zeta=\zeta_c$ (among a finite number of solutions)  to the equation above. In particular, $\zeta$ and $\omega$ can be thought of as homogeneous functions of degree zero in the variable $(r,s)$.

We define
\begin{eqnarray}
\label{def:a(z,theta)} a(\zeta,\theta)&:=&\frac{-G(\zeta e^{i\theta},g(\zeta e^{i\theta}))}{g(\zeta e^{i\theta})\,H_w(\zeta e^{i\theta},g(\zeta e^{i\theta}))}\,,\\
\label{def:f(z,theta)} f(\zeta,\theta)&:=&\ln\left\{\frac{g(\zeta e^{i\theta})}{g(\zeta)}\right\}-i\theta\frac{\zeta g'(\zeta)}{g(\zeta)}\,,
\end{eqnarray}
which are analytic for all $\theta$ sufficiently small and $\zeta$ sufficiently close to $\zeta_c$.

Our main result is the following one.

\begin{theorem}\label{thm:main}
Let $G(z,w)$, $H(z,w)$, $K$, $(\zeta_c,\omega_c)$, etc. be as above. Define $F(z,w):=G(z,w)/H(z,w)$. If there are nonnegative integers $p\le q$ such that $a(\zeta,\theta)$ has a $p$-to-$q$ change of degree about $\theta=0$ as $\zeta\to\zeta_c$, while $f(\zeta,\theta)$ vanishes to constant degree $n$ about $\theta=0$ as $\zeta\to\zeta_c$, then there is a constant $C>0$ and functions $A_k(\zeta)$ and $B_k(\zeta;s)$, with $p\le k\le q$, analytic in $\zeta$ near $\zeta=\zeta_c$, such that 
\begin{eqnarray}
\label{thm:Ak(z0)=0} A_k(\zeta_c)&=&0,\, p\le k< q\,,\\
\label{thm:Aq(z0)ne0} A_q(\zeta_c)&\ne&0\,,
\end{eqnarray}
and
\begin{equation}
\label{thm:[zrws]F}
[z^rw^s]F=\frac{\zeta^{-r}\omega^{-s}}{2\pi}\left\{\sum_{k=p}^q A_k(\zeta)\cdot B_k(\zeta;s)+O(e^{-s\cdot C})\right\}\,,
\end{equation}
uniformly for all $(r,s)\in\Lambda$ such that $r/s$ is sufficiently close to $d(\zeta_c,\omega_c)$. Furthermore, there are coefficients $c_k(\zeta;j)$, with $j\ge k$, which are analytic in $\zeta$ near $\zeta=\zeta_c$ such that each coefficient $B_k$ above admits an asymptotic expansion of the form
\begin{equation}\label{thm:Bk(z;s)}
B_k(\zeta;s)\approx \sum_{j=k}^\infty c_k(\zeta;j)\cdot\left(1+(-1)^j\cdot D(j,n)\right)\cdot\frac{1}{n}\Gamma\left(\frac{j+1}{n}\right)\cdot s^{-(j+1)/n},
\end{equation}
as $s\to\infty$, uniformly for all $(r,s)\in\Lambda$ such that $r/s$ is sufficiently close to $d(\zeta_c,\omega_c)$, where we have defined
\begin{equation}\label{thm:D(j,n)}
D(j,n):=\left\{\begin{array}{lcll}
1 & , & n\hbox{ even } &;\\
\exp\left(-\frac{i\pi(j+1)}{n}\cdot\hbox{sign}\{i\cdot[\theta^n]f(\zeta_c,\theta)\}\right) & , & n\hbox{ odd } &.
\end{array}\right.
\end{equation}
\end{theorem}

\noindent{\bf Remark 1:} The analytic coefficients $A_k$ in (\ref{thm:Ak(z0)=0}) and (\ref{thm:Aq(z0)ne0}) together with an auxiliary function $\alpha=\alpha(\zeta,\theta)$ are the unique analytic solutions (near $\zeta=\zeta_c$ and $\theta=0$) to the system of equations
\begin{eqnarray*}
\int_0^\theta a(\zeta,w)dw&=&\sum\limits_{k=p}^q\frac{A_k(\zeta)}{k+1}\alpha^{k+1}\,,\\
A_k(\zeta_c)&=&0\,,p\le k<q\,,\\
A_q(\zeta_c)&\ne&0\,,\\
\alpha&=&\alpha(\zeta,\theta)=\theta+\ldots
\end{eqnarray*}
In particular, it follows that
\begin{eqnarray}
\label{rmk:Ap(z)} A_p(\zeta)&=&[\theta^p]a(\zeta,\theta)\,,\\
\label{rmk:Aq(z0)} A_q(\zeta_c)&=&[\theta^q]a(\zeta_c,\theta)\,.
\end{eqnarray}
\noindent{\bf Remark 2:} In (\ref{thm:Bk(z;s)}) one has that
\begin{eqnarray}
\label{rmk:ck(z;s)} c_k(\zeta;k)&=&\left([\theta^n]f(\zeta,\theta)\right)^{-(k+1)/n}\,.
\end{eqnarray}
More generally, the coefficients $c_k(\zeta;j)$ are characterized by the identity $c_k(\zeta;j)=[\beta^j]\alpha^k\frac{\partial\alpha}{\partial\beta}$ where, for all $\zeta$ sufficiently close to $\zeta_c$, the variables $\alpha$ and $\beta$ are related to each other through the variable $\theta$ via the relations
\begin{eqnarray*}
\alpha &=& \alpha(\zeta,\theta),\\
\beta &=&\alpha\cdot\left([\theta^n]f(\zeta,\theta)\right)^{1/n}\cdot\left(1+\frac{f(\zeta,\theta)-\left([\theta^n]f(\zeta,\theta)\right)\,\alpha^n}{\left([\theta^n]f(\zeta,\theta)\right)\,\alpha^n}\right)^{1/n}.
\end{eqnarray*}

Theorem 1 is essentially equivalent to Theorem 3.3 in~\cite{PW01} when the amplitude term $a(\zeta,\theta)$ in~(\ref{def:a(z,theta)}) vanishes to constant degree in the variable $\theta$ about $\theta=0$ (the case $p=q$). The first application we show is concerned with precisely this case. The generating function in the following example is analyzed in~\cite{PW01}. However, here we perform a similar analysis but from the perspective of Theorem \ref{thm:main}.\\

\noindent{\bf Application 1. (Lattice Paths.)} The Delannoy numbers (see \cite{Stanley99}, pp. 185) are the coefficients $f_{r,s}$ that count the number of paths in the $\Z\times\Z$-lattice that join $(0,0)$ with $(r,s)$ with steps of the form $(0,1)$, $(1,1)$ and $(1,0)$. With the understanding that $f_{0,0}=1$ and $f_{r,s}=0$ whenever $r<0$ or $s<0$, it follows that $f_{r,s}=f_{r-1,s}+f_{r-1,s-1}+f_{r,s-1}$, for all integers $r,s\ge0$ except when $(r,s)=(0,0)$. Using this recursion it is almost direct to see that
\[F(z,w):=\sum_{r,s\ge0} f_{r,s}\,z^r\,w^s=\frac{1}{1-z-w-zw}\,.\]

The strictly minimal simple zeros of the denominator of $F$ are all of the form $(\zeta,\omega)$, with $\zeta\in(0,1)$ and $\omega=g(\zeta):=\frac{1-\zeta}{1+\zeta}$. Furthermore, one finds that
\[(r,s)\in\hbox{dir}(\zeta,\omega)\iff \frac{r}{s}=\frac{2\zeta}{1-\zeta^2}\,.\]
This allows an asymptotic analysis for $[z^rw^s]F$ as $(r,s)\to\infty$, uniformly for $(r,s)$ in any cone of the form $\Lambda=\{(r,s):d_1\le r/s\le d_2\}$, with $d_1>0$ and $d_2>0$ arbitrary constants. On the other hand, as shown in \cite{PW01}, one finds for $(r,s)\in\Lambda$ that
\[(r,s)\in\hbox{dir}(\zeta,\omega)\iff \zeta=\frac{\sqrt{r^2+s^2}-s}{r}\,,\omega=\frac{\sqrt{r^2+s^2}-r}{s}\,.\]

Using definitions (\ref{def:a(z,theta)}) and (\ref{def:f(z,theta)}) it follows that
\begin{eqnarray*}
a(\zeta,\theta)&=&\frac{1}{1-\zeta e^{i\theta}}\\
&=&\frac{1}{1-\zeta}+\frac{i\zeta}{(1-\zeta)^2}\theta+\ldots\\
f(\zeta,\theta)&=&\ln\left\{\frac{(1-\zeta e^{i\theta})(1+\zeta)}{(1+\zeta e^{i\theta})(1-\zeta)}\right\}+\frac{2i\zeta}{1-\zeta^2}\theta\\
&=&\frac{\zeta(1+\zeta^2)}{(1-\zeta^2)^2}\theta^2+\frac{i\zeta(1+6\zeta^2+\zeta^4)}{3(1-\zeta^2)^3}\theta^3+\ldots
\end{eqnarray*}
Since $a(\zeta,\theta)$ and $f(\zeta,\theta)$ respectively vanish to constant degree $0$ and $2$ at $\theta=0$, for all $\zeta\in(0,1)$, Theorem \ref{thm:main} implies that there is a constant $c>0$ and coefficients $B(r,s)$ such that
\begin{eqnarray*}
[z^rw^s]F&=&\frac{1}{2\pi}\left(\frac{\sqrt{r^2+s^2}-s}{r}\right)^{-r}\left(\frac{\sqrt{r^2+s^2}-r}{s}\right)^{-s}\cdot\left\{\frac{r\cdot B(r,s) }{r+s-\sqrt{r^2+s^2}}+O(e^{-s\cdot c})\right\}\,,\\
B(r,s)&=&2\sqrt{\pi}\frac{s}{r}\left(\frac{\sqrt{r^2+s^2}-s}{r}+\frac{r}{\sqrt{r^2+s^2}-s}\right)^{-1/2}\cdot s^{-1/2}+O(s^{-3/2})\,,
\end{eqnarray*}
uniformly for $(r,s)\in\Lambda$ as $(r,s)\to\infty$. In particular, it follows that
\[[z^rw^s]F\sim\left(\frac{\sqrt{r^2+s^2}-s}{r}\right)^{-r}\left(\frac{\sqrt{r^2+s^2}-r}{s}\right)^{-s}\cdot\sqrt{\frac{rs}{2\pi(r+s-\sqrt{r^2+s^2})^2\sqrt{r^2+s^2}}}\,,\]
whenever $(r,s)\to\infty$ at a comparable rate.$\Box$\\

Although the computations in Theorem~\ref{thm:main} can be involved, it gives a precise and unified understanding of the elements that are important to take into consideration when analyzing the asymptotic behavior of the coefficients of meromorphic functions in two variables. Furthermore, the calculations greatly simplify in situations where the coefficients $A_k(\zeta)$ are easily available. This is the main point of the following result which is a direct consequence of Remark 1.

\begin{corollary}\label{cor:q=p+1}
Under the hypothesis of Theorem \ref{thm:main} but for the special case in which $q=p+1$, if for all $\zeta$ sufficiently close to $\zeta_c$, $\theta(\zeta)$ is the only non-trivial solution of the equation: $a(\zeta,\theta)=0$, with $\theta$ in some open neighborhood of $\theta=0$, then
\begin{eqnarray}
\label{cor:ide:Ap(z)} A_p(\zeta)&=&[\theta^p]a(\zeta,\theta)\,,\\
\label{cor:ide:Ap+1(z)} A_{p+1}(\zeta)&=&\left(\frac{(-1)^{p+1}}{(p+1)(p+2)}\cdot\{[\theta^p]a(\zeta,\theta)\}^{p+2}\cdot\left\{\int_0^{\theta(\zeta)}a(\zeta,\xi)d\xi\right\}^{-1}\right)^{1/(p+1)}\,,
\end{eqnarray}
where the branch of the $(p+1)$-root above is to be selected so as to have $\lim\limits_{\zeta\to \zeta_c}A_{p+1}(\zeta)=[\theta^{p+1}]a(\zeta_c,\theta)$.
\end{corollary}

\noindent{\bf Application 2. (Lagrange Inversion Formula.)} If $t(x)$ is an analytic function of $x$ near $x=0$ such that $t(x)=x\cdot u(t(x))$, for a certain analytic function $u(x)$ with $u(0)\ne0$, then $[x^r]t(x)=[x^{r-1}](u(x))^r/r$ (see Section 5.4 in~\cite{Stanley99}). More generally, many problems related to the Lagrange Inversion Formula naturally lead to study the asymptotic behavior of coefficients of the form $[x^r](u(x))^sv(x)$, as $(r,s)\to\infty$ (see~\cite{Drm94} and~\cite{BFSS00}). These coefficients are related to those of a bivariate generating function via the identity
\begin{equation}\label{exa:[xr]u(x)sv(x)}
[x^r](u(x))^sv(x)=[z^rw^s]\frac{v(z)}{1-wu(z)}\,.
\end{equation}
(See final remark on Section 2 in \cite{BFSS01} and remark 5.22 in \cite{Lla03} for the uses of multivariate generating functions in problems associated to the Lagrange Inversion Formula. See \cite{Wilson05} for a discussion in the context of Riordan arrays.)

In what follows we assume that the radius of convergence of $v(z)$ is greater or equal to that of $u(z)$. In the context of Theorem \ref{thm:main}, a point of the form $(\zeta,1/u(\zeta))$ is a strictly minimal simple zero of the denominator in the right-hand side of (\ref{exa:[xr]u(x)sv(x)}) provided that $\zeta\cdot u(\zeta)\ne0$ and that $|u(x)|$ is maximized on the circumference $|x|=|\zeta|$ solely at $x=\zeta$. We emphasize that this condition is easily satisfied for $\zeta>0$ and within the radius of convergence of $u(z)$ whenever $u(z)$ is aperiodic and has non-negative Taylor coefficients. Asymptotic formulae for the coefficients in (\ref{exa:[xr]u(x)sv(x)}) are then available along the directions in the $(r,s)$-lattice where $r/s=\zeta u'(\zeta)/u(\zeta)$. Furthermore, if $K$ is a compact set of strictly minimal simple zeros and $(\zeta_c,1/u(\zeta_c))$ is an interior point of $K$, then Theorem \ref{thm:main} can be used to provide asymptotic formulae in an open cone of directions in the $(r,s)$-lattice containing the line $r/s=\zeta_c u'(\zeta_c)/u(\zeta_c)$, provided that there are non-negative integers $p\le q$ and $n$ such that
\begin{eqnarray}
\label{exa:Lagrange a(z,theta)} a(\zeta,\theta)&:=& v(\zeta e^{i\theta})\,,\\
\label{exa:Lagrange f(z,theta)} f(\zeta,\theta)&:=& \ln\left\{\frac{u(\zeta)}{u(\zeta e^{i\theta})}\right\}+i\theta\frac{\zeta u'(\zeta)}{u(\zeta)}\,,
\end{eqnarray}
respectively have a $p$-to-$q$ and $n$-to-$n$ change of degree about $\theta=0$ as $\zeta\to\zeta_c$.\\

\noindent To fix ideas consider the case in which $u(x):=(1-x)^{-1}$ and $v(x):=(1-2x)$. Then every point of the form $(\zeta,1-\zeta)$, with $\zeta\in(0,1)$, is a strictly minimal simple zero of the denominator in the right-hand side of (\ref{exa:[xr]u(x)sv(x)}). Furthermore, $(r,s)\in\dir(\zeta,1-\zeta)$ if and only if $r/s=\zeta/(1-\zeta)$; in particular,
\[(r,s)\in\dir(\zeta,1-\zeta)\iff \zeta=\frac{r}{r+s}\,.\]
This motivates us to define $\zeta(r,s):=r/(r+s)$, for all $(r,s)$ such that $r\cdot s>0$.

Observe that back in (\ref{exa:Lagrange a(z,theta)}) and (\ref{exa:Lagrange f(z,theta)}) one finds that
\begin{eqnarray*}
a(\zeta,\theta)&=&(1-2\zeta)-2i\zeta\theta+\ldots\\
f(\zeta,\theta)&=&\frac{\zeta}{2(1-\zeta)^2}\theta^2+\ldots
\end{eqnarray*}
While $f(\zeta,\theta)$ vanishes to constant degree 2 about $\theta=0$, for all $\zeta\in(0,1)$, $a(\zeta,\theta)$ has a $0$-to-$1$ change of degree about $\theta=0$, as $\zeta\to1/2$. As a result, using Theorem 1, we can determine the asymptotic behavior of $[x^r](1-x)^{-s}(1-2x)$ as $(r,s)\to\infty$ so long as $r$ and $s$ grow at a comparable rate. 

Theorem 1 implies almost immediately that
\begin{equation}\label{exa:away diagonal}
[x^r](1-x)^{-s}(1-2x)=\frac{1}{\sqrt{2\pi}}\left(\frac{r}{r+s}\right)^{-r}\left(\frac{s}{r+s}\right)^{-s}\cdot\left\{\left(1-\frac{r}{s}\right)\left(1+\frac{s}{r}\right)^{1/2}s^{-1/2}+O(s^{-1})\right\}\,,
\end{equation}
as $(r,s)\to\infty$, uniformly for $r/s$ restricted to a compact subset of $(0,1)\cup(1,\infty)$.

On the other hand, Corollary~\ref{cor:q=p+1} implies that there is an $\epsilon>0$ such that
\[[x^r](1-x)^{-s}(1-2x)=\frac{1}{\sqrt{2\pi}}\left(\frac{r}{r+s}\right)^{-r}\left(\frac{s}{r+s}\right)^{-s}
\cdot\left\{A_0\left(\frac{r}{r+s}\right)\cdot B_0(r,s)+A_1\left(\frac{r}{r+s}\right)\cdot B_1(r,s)\right\}\,,\]
as $(r,s)\to\infty$, uniformly for $(1-\epsilon)\le r/s\le(1+\epsilon)$, where
\begin{eqnarray*}
A_0(\zeta)&:=& 1-2\zeta\,,\\
B_0(r,s) &=& \alpha_0\left(\frac{r}{r+s}\right)s^{-1/2}+O(s^{-3/2})\,,\\
\alpha_0(\zeta)&:=&\frac{1-\zeta}{\sqrt{\zeta}}\,,\\
A_1(\zeta)&:=&\frac{i(1-2\zeta)^2}{2(1-2\zeta+\ln(2\zeta))}\,,\\
B_1(r,s) &=& \alpha_1\left(\frac{r}{r+s}\right)s^{-3/2}+O(s^{-5/2})\,,\\
\alpha_1(\zeta)&:=&-i(1-\zeta)^2\cdot\frac{5-9\zeta-12\zeta^2+20\zeta^3+2(1+5\zeta-8\zeta^2)\ln(2\zeta)}{4\zeta\sqrt{\zeta}(1-2\zeta)(1-2\zeta+\ln(2\zeta))}\,.
\end{eqnarray*}
Observe that Theorem 1 asserts that $A_1(\zeta)$ and $\alpha_1(\zeta)$ are analytic about any $\zeta\in(0,1)$. The apparent singularity of $A_1(\zeta)$ at $\zeta=1/2$ is not such because its denominator vanishes to degree 2 about $\zeta=1/2$. On the other hand, the numerator and denominator of $\alpha_1(\zeta)$ vanish to degree 3 about $\zeta=1/2$. Indeed, the first few terms of the Taylor series of $A_1(\zeta)$ and $\alpha_1(\zeta)$ about $\zeta=1/2$ are found to be
\begin{eqnarray*}
A_1(\zeta)&=&-i-\frac{4i}{3}\left(\zeta-\frac{1}{2}\right)+\frac{2i}{9}\left(\zeta-\frac{1}{2}\right)^2+\ldots\\
\alpha_1(\zeta)&=&-\frac{i\sqrt{2}}{4}+\frac{31i\sqrt{2}}{24}\left(\zeta-\frac{1}{2}\right)-\frac{503i\sqrt{2}}{180}\left(\zeta-\frac{1}{2}\right)^2+\ldots
\end{eqnarray*}

Since $A_0(r/(r+s))=0$ whenever $r=s$, the above expansion for $[x^r](1-x)^{-s}(1-2x)$ implies that
\begin{equation}\label{exa:along diagonal}
[x^r](1-x)^{-s}(1-2x)=-\frac{4^{(s-1)}}{\sqrt{\pi}}\left\{s^{-3/2}+O(s^{-5/2})\right\}\,,
\end{equation}
as $(r,s)\to\infty$ with $r=s$. This corresponds to the asymptotic expansion one would obtain after using Stirling's formula to find the leading asymptotic order of the factorial terms in the identity
\[[x^r](1-x)^{-r}(1-2x)=-\frac{(2r-2)!}{r((r-1)!)^2}\,.\]

Formulae (\ref{exa:away diagonal}) and (\ref{exa:along diagonal}) characterize the asymptotic behavior of the coefficients $[x^r](1-x)^{-s}(1-2x)$ as $(r,s)\to\infty$ along the diagonal line $r=s$ or along directions completely away from it. More explicit asymptotic formulae for these coefficients, as $r/s\to1$, can be obtained looking at 
the Taylor coefficients of the functions $A_0(\zeta)\cdot\alpha_0(\zeta)$ and $A_1(\zeta)\cdot\alpha_1(\zeta)$ about $\zeta=1/2$. Indeed, it follows for all constant $\delta>0$ that
\begin{equation}\label{exa:near diagonal}
[x^r](1-x)^{-s}(1-2x)=\frac{-1}{2\sqrt{\pi}}\left(\frac{r}{r+s}\right)^{-r}\left(\frac{s}{r+s}\right)^{-s}
\cdot\left\{\frac{r-s}{r+s}\cdot s^{-1/2}+\frac{s^{-3/2}}{2}+O(s^{-5/2})\right\}\,,
\end{equation}
as $(r,s)\to\infty$, uniformly for $(r,s)$ in the region $1-\delta/s\le r/s\le 1+\delta/s$.

If $r/s$ approaches $1$ from above then cancellation between the first two terms in the curly bracket above is ruled out. As a result, if $r/s=1+|O(s^{-1})|$ then 
\begin{equation}\label{exa:above diagonal}
[x^r](1-x)^{-s}(1-2x)\sim\frac{-1}{2\sqrt{\pi}}\left(\frac{r}{r+s}\right)^{-r}\left(\frac{s}{r+s}\right)^{-s}
\cdot\left\{\frac{r/s-1}{r/s+1}\cdot s+\frac{1}{2}\right\}\cdot s^{-3/2}\,.
\end{equation}
This means that in the $(r,s)$-lattice a bandwidth of size $s^{-1}$ from above the line 
$r=s$ is what separates the behavior of $[x^r](1-x)^{-s}(1-2x)$ as prescribed in (\ref{exa:away diagonal}) from the one in (\ref{exa:along diagonal}).

On the other hand, if $r/s$ approaches $1$ from below then a cascade effect of cancellation in (\ref{exa:near diagonal}) may reduce the size of $[x^r](1-x)^{-s}(1-2x)$ to arbitrarily small orders. Refined estimates in this case depend on the precise rate of convergence of $r/s$ toward 1. To amplify on this consider coefficients $\alpha>0$, $\beta\ge1$, $\gamma\ne0$ and $\delta>0$ and suppose that
\[\frac{r}{s}=1-\alpha s^{-\beta}+\gamma s^{-(\beta+\delta)}+o(s^{-(\beta+\delta)})\,.\]
In particular, $(r-s)/(r+s)=-\alpha s^{-\beta}(1+\alpha s^{-\beta}/2)/2+\gamma s^{-(\beta+\delta)}/2+o(s^{-2\beta}+s^{-(\beta+\delta)})$. Using this in (\ref{exa:near diagonal}) we obtain that
\[[x^r](1-x)^{-s}(1-2x)\sim\frac{-1}{4\sqrt{\pi}}\left(\frac{r}{r+s}\right)^{-r}\left(\frac{s}{r+s}\right)^{-s}
\cdot \left\{\begin{array}{ll}
(1-\alpha)s^{-3/2} &\hbox{, if }\alpha\ne1\hbox{ and }\beta=1;\\
\gamma s^{-(3/2+\delta)} &\hbox{, if }\alpha=1,\beta=1\hbox{ and }0<\delta<1;\\
s^{-3/2}&\hbox{, if }\beta>1.
\end{array}
\right.\]
As a result and unlike the asymptotic description in (\ref{exa:above diagonal}), we see that if $(r,s)\to\infty$ with $r/s\uparrow1$ then there is no well-defined bandwidth that separates the asymptotic behavior of $[x^r](1-x)^{-s}(1-2x)$ as prescribed in (\ref{exa:away diagonal}) from the one in (\ref{exa:along diagonal}). Furthermore, if $\alpha=\beta=1$ and $0<\delta<1$ then $[x^r](1-x)^{-s}(1-2x)$ is of an asymptotic order smaller than anyone observed as $(r,s)\to\infty$ along any diagonal line in the $(r,s)$-lattice. This finding is consistent with the identity
\[[x^r](1-x)^{-s}(1-2x)=\frac{(s-r-1)\cdot(r+s-2)!}{r!\cdot(s-1)!}\,,\]
from which we see that $[x^r](1-x)^s(1-2x)=0$ whenever $r/s=1-s^{-1}$.$\Box$\\

\noindent{\bf Remark 3:} The determination of the coefficients $A_k(\zeta)$ in Theorem~\ref{thm:main} becomes more difficult the bigger is the change of degree of the amplitude term $a(\zeta,\theta)$ in (\ref{def:a(z,theta)}). However, the linear dependence between the asymptotic expansion of the coefficients of $F$ and of $a(\zeta,\theta)$ can be exploited to overcome this problem.
Indeed, if $a(\zeta,\theta)$ has a $p$-to-$q$ change of degree in the variable $\theta$, with $p<q$, then one can rewrite $a(\zeta,\theta)=a_0(\zeta,\theta)+a_1(\zeta,\theta)$, where $a_0(\zeta,\theta)$ is a polynomial in the variable $\theta$ (of degree less than $q$) and $a_1(\zeta,\theta)$ vanishes regardless of $\zeta$ to constant degree $q$ in $\theta$. Theorem 1 can now be used to obtain an asymptotic expansion for each of the terms in $a_0(\zeta,\theta)$ as well as for $a_1(\zeta,\theta)$. Combining these linearly, one obtains an asymptotic expansion for $[z^rw^s]F$ that resembles the one in~(\ref{thm:[zrws]F}).

\section{Proof of Main Results}

\subsection{Associating a parameter-varying integral}

In this section we show some preliminary results that are required to prove Theorem~\ref{thm:main}. We assume that there are functions $G(z,w)$ and $H(z,w)$ analytic in an open polydisk $D$ centered at $(0,0)$ on which $F(z,w)$, the generating function associated to the coefficients $(f_{r,s})$, satisfies the identity $F(z,w)=G(z,w)/H(z,w)$. In addition, we assume as given a compact set $K\subset D$ of strictly minimal simple zeros of $H$ containing a particular point $(\zeta_c,\omega_c)$. It is assumed that $H_w(\zeta_c,\omega_c)\ne0$. In particular, the Implicit Function Theorem implies that $(\zeta_c,\omega_c)$ has an open neighborhood of the form $Z\times W\subset D$ and there is an analytic map $g:Z\to W$ such that for all $(z,w)\in Z\times W$,
$H(z,w)=0\hbox{ if and only if }w=g(z)$. Without loss of generality we may assume that $0\notin W$.

We now adopt the following notation. For all $0<\epsilon<\pi/2$, the notation $|\arg\{z\}|\le\epsilon$ signifies that $z=|z|e^{i\theta}$, for some $\theta\in[-\epsilon,\epsilon]$. Accordingly, the notation $|\arg\{z\}|\ge\epsilon$ is used to mean that $z=|z|e^{i\theta}$, for some $\theta\in[\epsilon,\pi]\cup[-\epsilon,-\pi]$.

\begin{lemma}\label{lem:H zeroe free 1} For all $\epsilon_1>0$ sufficiently small
there is a $\delta_1>0$ such that for all $(\zeta,\omega)\in K$, $H$ is
zero-free on the set
$\{(z,w):|z|=|\zeta|,|\arg(z/\zeta)|\ge\epsilon_1,|w|\le(1+\delta_1)|\omega|\}$.
\end{lemma}

\noindent{\bf Proof:} Without loss of generality, assume that
$0<\epsilon_1<\pi/2$. If $K$ consisted of only one point the lemma
would follow directly from the continuity of $H$ together with the
strict minimality of its only element. More generally, define for
each $(\zeta,\omega)\in K$ the quantity $\delta_1(\zeta,\omega)$ to be the supremum
of those $\delta>0$ such that $H$ is zero-free on the set
$\{(z,w):|z|=|\zeta|,|\arg(z/\zeta)|\ge\epsilon_1,|w|\le(1+\delta)|\omega|\}$.
To prove the lemma it is enough to show that $\inf\{\delta_1(\zeta,\omega):(\zeta,\omega)\in K\}>0$. We prove this by contradiction. Assuming otherwise there would be a sequence of points $(\zeta_j,\omega_j)\in K$ such that
$\delta_1(\zeta_j,\omega_j)\to0$, as $j\to\infty$. In particular, for
all $j$ sufficiently large, there would be a $(z_j,w_j)$ 
such that $|z_j|=|\zeta_j|$, $|\arg\{z_j/\zeta_j\}|\ge\epsilon_1$,
$|w_j|=(1+\delta_1(\zeta_j,\omega_j))|\omega_j|$, and $H(z_j,w_j)=0$.
But, since $K$ is a compact set, there is no loss of generality in assuming that $(\zeta_j,\omega_j)\to(\zeta,\omega)\in K$ and $(z_j,w_j)\to(z,w)$, as $j\to\infty$. In particular, $|z|=|\zeta|$, $|w|=|\omega|$, $z\ne\zeta$, however, $H(z,w)=0$. This contradicts the fact that $(\zeta,\omega)$ is a strictly minimal zero of $H$ and therefore we conclude that $\inf\{\delta_1(\zeta,\omega):(\zeta,\omega)\in K\}>0$. This completes the proof of the lemma.$\Box$

\begin{lemma}\label{lem:H zeroe free 2}
For all $\epsilon_2>0$ sufficiently small there is a $\delta_2>0$ such
that all zeros of $H$ in the set $\{(z,w):|z-\zeta_c|<\epsilon_2,|w|\le(1+\delta_2)|g(z)|\}$ are of the form $w=g(z)$.
\end{lemma}

\noindent{\bf Proof:} The strict minimality of $(\zeta_c,\omega_c)$ together with the analyticity of $H$ imply that there is $\eta>0$ such that $w=\omega_c$ is the only zero of $H(\zeta_c,w)$ in the disk $\{w:|w|\le(1+\eta)|\omega_c|\}$ (see Theorem 10.18 in~\cite{Rudin87}). Without loss of generality we may assume that $\{w:|w-\omega_c|\le\eta|\omega_c|\}\subset W$. Observe that $H(\zeta_c,w)$ is zero-free on the set $\{w:\eta|\omega_c|\le|w-\omega_c|\hbox{ and }|w|\le(1+\eta)|\omega_c|\}$. Thus, since $H$ is uniformly continuous, 
it follows for all $\epsilon_2>0$ sufficiently small that $H$ is zero-free in the set $\{(z,w):|z-\zeta_c|\le\epsilon_2\,,\eta|\omega_c|\le|w-\omega_c|\hbox{ and }|w|\le(1+\eta)|\omega_c|\}$. In this case, the condition that $\{w:|w-\omega_c|\le\eta|\omega_c|\}\subset W$ implies that all zeros of $H$ in the polydisk $\{(z,w): |z-\zeta_c|\le\epsilon_2, |w|\le(1+\eta)|\omega_c|\}$ are of the form $w=g(z)$. The lemma follows after selecting $\epsilon_2>0$ small enough and defining $\delta_2>0$ so as to have
\[(1+\delta_2)=(1+\eta)\inf\limits_{z:|z-\zeta_c|\le\epsilon_2}\left|\frac{\omega_c}{g(z)}\right|>1\,.\]
The above inequality is always possible because $g(\zeta_c)=\omega_c$. This completes the proof of the lemma.$\Box$\\

The next result pretty much follows the lines of Lemma 4.1 in~\cite{PW01}. It is included in here for the sake of completeness.

\begin{lemma}\label{lem:[zrws]F as an integral}
For a sufficiently small choice of $\epsilon>0$ and for all $|\theta|\le\epsilon$ and $\zeta$ sufficiently close to $\zeta_c$, consider the functions $a(\zeta,\theta)$ and $f(\zeta,\theta)$ as defined in (\ref{def:a(z,theta)}) and (\ref{def:f(z,theta)}) respectively. Then $f(\zeta,0)=\frac{\partial f}{\partial\theta}(\zeta,0)=0$. In addition, for all $(\zeta,\omega)\in K$ sufficiently close to $(\zeta_c,\omega_c)$ and for all nonzero $\theta$ such that $-\epsilon\le\theta\le\epsilon$, $\Re\{f(\zeta,\theta)\}>0$. Furthermore, if
\begin{equation}
\label{lem:Sigma(z;s)} \Sigma(\zeta;s):=\int_{-\epsilon}^\epsilon e^{-s\cdot f(\zeta,\theta)} a(\zeta,\theta)\,d\theta\,,
\end{equation}
then there is a constant $c>0$ such that
\begin{equation}\label{lem:frs asymptotic formula}
[z^rw^s]F=\frac{\zeta^{-r}\omega^{-s}}{2\pi}\left\{\Sigma(\zeta;s)+O(e^{-sc})\right\}\,,
\end{equation}
uniformly for all $(r,s)\in dir(\zeta,\omega)$ and $(\zeta,\omega)\in K$ sufficiently close to $(\zeta_c,\omega_c)$.
\end{lemma}

\noindent{\bf Proof:} Let $\epsilon_2>0$ and $\delta_2>0$ be as in Lemma~\ref{lem:H zeroe free 2}. Consider $\epsilon_3>0$ be such that the functions $a(\zeta,\theta)$ and $f(\zeta,\theta)$ are analytic for $|\zeta-\zeta_c|\le\epsilon_3$ and $|\theta|\le\epsilon_3$. In addition, consider for $\epsilon_1>0$ the sets
\begin{eqnarray*}
K_c&:=&\{(\zeta,\omega)\in K: |\zeta-\zeta_c|\le\epsilon_1\}\,,\\
\gamma_1(\zeta) &:=& \{z:|z|=|\zeta|\hbox{ and }|\arg\{z/\zeta\}|\ge\epsilon_1\}\,,\\
\gamma_2(\zeta) &:=& \{z:|z|=|\zeta|\hbox{ and }|\arg\{z/\zeta\}|\le\epsilon_1\}\,.
\end{eqnarray*}
Select $\epsilon_1>0$ small enough so as to have $\gamma_2(\zeta)\subset\{z:|z-\zeta_c|\le\min\{\epsilon_2,\epsilon_3\}\}$, whenever $(\zeta,\omega)\in K_c$. Furthermore, chose $\epsilon_1>0$ sufficiently small so that the conclusion of Lemma \ref{lem:H zeroe free 1} applies with some $\delta_1>0$. Select $\delta$ so as to satisfy $0<\delta<\min\{\delta_1,\delta_2,1\}$. The strict minimality of $(\zeta,\omega)\in K_c$ implies that $H$ is zero-free on the polydisk $\{z:|z|\le|\zeta|\}\times\{w:|w|\le(1-\delta)|\omega|\}$. Cauchy's Formula \cite{Rudin87} then can be used to represent the coefficients of $F$ by the integrals
\begin{equation}\label{proof:lem:Cauchy formula 1}
[z^rw^s]F=\frac{1}{2\pi}\left\{\int_{z\in\gamma_1(\zeta)}+\int_{z\in\gamma_2(\zeta)}\right\}\frac{1}{z^r}\left(\frac{1}{2\pi i}\int_{|w|=(1-\delta)|\omega|}\frac{G(z,w)}{H(z,w)\cdot w^{s+1}}dw\right)\frac{dz}{iz}\,,
\end{equation}
where all contour integrals are in the standard counter-clockwise orientation. 

Lemma~\ref{lem:H zeroe free 1} implies  for all $(\zeta,\omega)\in K_c$ that $H$ is zero-free on the set $\gamma_1(\zeta)\times\{w:|w|\le(1+\delta)|\omega|\}$. As a result,
\begin{eqnarray*}
\left|\int_{z\in\gamma_1(\zeta)}\frac{1}{z^r}\int_{w=(1-\delta)|\omega|}\frac{G(z,w)}{H(z,w)\cdot w^{s+1}}dw\frac{dz}{iz}\right|
&=&\left|\int_{z\in\gamma_1(\zeta)}\frac{1}{z^r}\int_{w=(1+\delta)|\omega|}\frac{G(z,w)}{H(z,w)\cdot w^s}\frac{dw}{iw}\frac{dz}{iz}\right|\,,\\
&\le& (2\pi)^2|\zeta|^{-r}\{(1+\delta)|\omega|\}^{-s}\cdot\sup\limits_{\Gamma_1}|F|\,,
\end{eqnarray*}
where for convenience we have defined $\Gamma_1$ to be the set of all those points of the form $(z,w)$ such that there exists a $(\zeta,\omega)\in K_c$ such that $z\in\gamma_1(\zeta)$ and $|w|\le(1+\delta)|\omega|$. Since $\Gamma_1$ is compact and $H$ is zero-free over it, then $\sup\limits_{\Gamma_1}|F|$ must be finite. Back in (\ref{proof:lem:Cauchy formula 1}), this implies that
\begin{equation}\label{proof:lem:Cauchy formula 2}
[z^rw^s]F=\frac{1}{2\pi}\int_{z\in\gamma_2(\zeta)}\frac{1}{z^r}\left(\frac{1}{2\pi i}\int_{|w|=(1-\delta)|\omega|}\frac{G(z,w)}{H(z,w)\cdot w^{s+1}}dw\right)\frac{dz}{iz}+O(|\zeta|^{-r}|\omega|^{-s}(1+\delta)^{-s})\,,
\end{equation}
uniformly for all $r,s\ge0$ and all $(\zeta,\omega)\in K_c$. However, Lemma~\ref{lem:H zeroe free 2} implies that for each $(\zeta,\omega)\in K_c$ and $z\in\gamma_2(\zeta)$, $w=g(z)$ is the only singularity of the integrand above within the disk $\{w:|w|\le(1+\delta)|g(z)|\}$. The Residue Theorem in one variable~\cite{Rudin87} lets us conclude that
\[\frac{1}{2\pi i}\int_{|w|=(1-\delta)|\omega|}\frac{G(z,w)}{H(z,w)\cdot w^{s+1}}dw=\frac{-G(z,g(z))}{H_w(z,g(z))\cdot\{g(z)\}^{s+1}}+\frac{1}{2\pi i}\int_{|w|=(1+\delta)|g(z)|}\frac{G(z,w)}{H(z,w)\cdot w^{s+1}}dw\,.\]
But, observe that if $|w|=(1+\delta)|g(z)|$ and $z\in\gamma_2(\zeta)$ then the strict minimality of $(\zeta,\omega)\in K_c$ implies that $|g(z)|\ge|g(\zeta)|=|\omega|$. In particular,
\[\left|\frac{1}{2\pi}\int_{z\in\gamma_2(\zeta)}\frac{1}{z^r}\left(\frac{1}{2\pi i}\int_{|w|=(1+\delta)|g(z)|}\frac{G(z,w)}{H(z,w)\cdot w^{s+1}}dz\right)\frac{dz}{iz}\right|\le
|\zeta|^{-r}\{(1+\delta)|\omega|\}^{-s}\cdot\sup_{\Gamma_2}|F|\,,\]
where we have defined $\Gamma_2$ to be the set of points $(z,w)$ for which there exists a $(\zeta,\omega)\in K_c$ such that $z\in\gamma_2(\zeta)$ and $|w|=(1+\delta)|g(z)|$. Since $\Gamma_2$ is a compact set and $H$ is zero-free over it, from (\ref{proof:lem:Cauchy formula 2}) we can conclude that
\begin{equation}\label{proof:lem:[zrws]F formula 1}
[z^rw^s]F=\frac{1}{2\pi}\int_{z\in\gamma_2(\zeta)}\frac{1}{z^r}\frac{-G(z,g(z))}{H_w(z,g(z))\cdot\{g(z)\}^{s+1}}\frac{dz}{iz}+O(|\zeta|^{-r}|\omega|^{-s}(1+\delta)^{-s})\,,
\end{equation}
uniformly for all $r,s\ge0$ and all $(\zeta,\omega)\in K_c$.

The integral on the right-hand side in (\ref{proof:lem:[zrws]F formula 1}) can be parametrized using polar coordinates. Indeed, substituting $z=\zeta e^{i\theta}$, with $-\epsilon\le\theta\le\epsilon$, one 
obtains that
\[[z^rw^s]F=\frac{\zeta^{-r}\omega^{-s}}{2\pi}\int_{-\epsilon}^\epsilon e^{-s\cdot f(\theta;\zeta,r/s)}a(\zeta,\theta)d\theta+O(|\zeta|^{-r}|\omega|^{-s}(1+\delta)^{-s})\,,\]
uniformly for all $r,s\ge0$ and all $(\zeta,\omega)\in K_c$, where $a(\zeta,\theta)$ is defined as in (\ref{def:a(z,theta)}) and
$f(\theta;\zeta,\lambda):=\ln\left\{\frac{g(\zeta e^{i\theta})}{g(\zeta)}\right\}+i\lambda\theta$.
Observe that $f\left(0;\zeta,\frac{r}{s}\right)=0$ and
\begin{eqnarray*}
\frac{\partial f}{\partial\theta}\left(0;\zeta,\frac{r}{s}\right)&=&i\left(\frac{\zeta g'(\zeta)}{g(\zeta)}+\frac{r}{s}\right)\,,\\
&=&i\left(\frac{r}{s}-\frac{\zeta H_z(\zeta,\omega)}{\omega H_w(\zeta,\omega)}\right)\,.
\end{eqnarray*}
In particular we see that $\frac{\partial f}{\partial\theta}\left(0;\zeta,\frac{r}{s}\right)=0$, for all $(r,s)\in\hbox{dir}(\zeta,\omega)$.  Furthermore, the strict minimality of $(\zeta,\omega)\in K_c$ implies that $|g(\zeta e^{i\theta})|>|g(\zeta)|$, for all nonzero $\theta$ such that $-\epsilon\le\theta\le\epsilon$ and, as a result, $\Re\{f(\zeta;\theta)\}>0$ for all such $\theta$. The lemma follows by noticing that whenever $(r,s)\in\dir(\zeta,\omega)$ then $f\left(0;\zeta,r/s\right)=f(\zeta;\theta)$, with $f(\zeta;\theta)$ as defined in (\ref{def:f(z,theta)}).$\Box$

\subsection{Polynomial Canonical Representations.} A result of Levinson~\cite{Lev60III} implies that if a function $H(u,v)$ is analytic in a neighborhood of the origin in $\CC^2$ and its Hartog's series vanishes to degree $q\ge1$ in the variable $v$ about the origin, then $H$ admits a near $(0,0)$ a representation of the form
\begin{equation}\label{ide:Levison Lladser representation}
H(u,v)=\sum_{j=0}^q H_j(u)w^j\,.
\end{equation}
Above the coefficient functions $H_j$ are analytic near the origin and such that $H_j(0)=0$, for all $0\le j<q$, however, $H_q(0)\ne0$. In addition, $w=w(u,v)$ is a certain analytic function near the origin such that $w(u,0)=0$ and $\frac{\partial w}{\partial v}(u,0)=1$. In \cite{Lla03} it is proved using one complex variable methods that this representation is indeed unique. The following more precise representation will be more suitable to prove our main result.

\begin{lemma}\label{lem:p-to-q terms}
Let $0\le p\le q$ with $q\ge1$ be nonnegative integers. Suppose that $H(u,v)$ is analytic in a neighborhood of the origin and it has a $p$-to-$q$ change of degree about $v=0$ as $u\to0$. Then, $H$ admits near the origin a unique representation of the form
\begin{equation}\label{lem:ide:p-to-q terms}
H(u,v)=\sum_{k=p}^q H_k(u)\cdot w^k\,,
\end{equation}
where $H_k(0)=0$, for $p\le k<q$, $H_q(0)\ne0$, and $w=w(u,v)$ is such
that $w(u,0)=0$ and $\frac{\partial w}{\partial v}(u,0)=1$. Furthermore,
\begin{equation}\label{lem:ide:pth term in Levinson Lladser form}
H_p(u)=\frac{1}{p!}\frac{\partial^pH}{\partial v^p}(u,0)\,.
\end{equation}
\end{lemma}

\noindent{\bf Proof:} The uniqueness of the representation in (\ref{lem:ide:p-to-q terms}) is immediate from the uniqueness of the representation in (\ref{ide:Levison Lladser representation}). Suppose that the representation in (\ref{ide:Levison Lladser representation}) applies for all $(u,v)$ in an open neighborhood of the polydisk $\{(u,v):|u|\le\epsilon\hbox{ and }|v|\le\epsilon\}$, for some $\epsilon>0$.  Since $H(u,v)$ has a $p$-to-$q$ change of degree about $v=0$ as $u\to0$, $H$ has a Hartog's series of the form
\[H(u,v)=\sum\limits_{k=p}^\infty h_k(u)\,v^k\,,\]
where the coefficients $h_k$ are analytic for $|u|\le\epsilon$, $h_p(u)$ is not identically zero in any neighborhood of $u=0$, and $h_q(0)\ne0$.

Consider the map $\Phi(u,v)=(u,w(u,v))$.  The conditions imposed over $w$ in (\ref{ide:Levison Lladser representation}) imply that the Jacobian matrix $\frac{\partial\Phi}{\partial(u,v)}(0,0)$ is triangular with all entries
equal to 1 along the diagonal. Since $\Phi(0,0)=(0,0)$, the Inverse Mapping Theorem
lets us assume without loss of generality that $\Phi$ is holomorphic and $1$-to-$1$ over the polydisk $\{(u,v):|u|\le\epsilon,|v|\le\epsilon\}$. In particular, for all $u$ such that $|u|\le\epsilon$, $w(u,\cdot)$ is 1-to-1 for $|v|\le\epsilon$. Furthermore, since $w(u,0)=0$, the Open Mapping Theorem implies that there are $\rho_1,\rho_2>0$ such that $\{w:|w|\le\rho_1\}\subset w(u,\{v:|v|<\epsilon\})$ and  the pre-image of
$\{v:|v|<\rho_1\}$ under $w(u,\cdot)$ contains the disk $\{v:|v|\le\rho_2\}$. As a result, using Cauchy's Formula in (\ref{ide:Levison Lladser representation}) and then the substitution $w=w(u,v)$, it follows for all $0\le j\le q$, that
\begin{eqnarray*}
H_j(u)&=&\frac{1}{2\pi i}\int_{|w|=\rho_1}\frac{1}{w^{j+1}}\left(\sum_{k=0}^q H_k(u)\cdot w^k\right)\,dw\,,\\
&=&\frac{1}{2\pi i}\int_{|v|=\rho_2}\frac{H(u,v)}{\{w(u,v)\}^{j+1}}\,\frac{\partial
w}{\partial v}(u,v)\,dv\,,\\
&=&\frac{1}{2\pi i}\sum_{k=p}^\infty
h_k(u)\cdot\int_{|v|=\rho_2}\frac{v^k}{\{w(u,v)\}^{j+1}}\,\frac{\partial w
}{\partial v}(u,v)\,dv\,,
\end{eqnarray*}
where for the last identity we have used that the Hartog's series of $H$ converges uniformly over compact subsets of $\{(u,v):|u|\le\epsilon\hbox{ and }|v|\le\epsilon\}$. However, observe that the conditions imposed over $w$ in (\ref{ide:Levison Lladser representation}) imply that, for all $j<p\le k$, the function $\frac{v^k}{\{w(u,v)\}^{j+1}}\,\frac{\partial w}{\partial v}(u,v)$ is
analytic in $v$ in an open neighborhood of $\{v:|v|\le\rho_2\}$. Consequently, for $j<p$, all the terms in the above summation vanish and therefore $H_j(u)=0$. This shows (\ref{lem:ide:p-to-q terms}). Furthermore, if $j=p$ then the Residue Theorem implies that
\begin{eqnarray*}
H_p(u)&=&\frac{h_p(u)}{2\pi i}\cdot\int_{|v|=\rho_2}\frac{v^p}{\{w(u,v)\}^{p+1}}\,\frac{\partial w}{\partial v}(u,v)\,dv\,,\\
&=&h_p(u)\cdot\hbox{Res}\left(\frac{v^p}{\{w(u,v)\}^{p+1}}\,\frac{\partial w}{\partial v}(u,v);v=0\right)\,,\\
&=&h_p(u)\,.
\end{eqnarray*}
This shows (\ref{lem:ide:pth term in Levinson Lladser form}) and completes the proof of the lemma.$\Box$

\subsection{Asymptotic Analysis.}\label{sec:asymp analy}
In this section we prove Theorem \ref{thm:main}. This is accomplished by analyzing the asymptotic behavior of the integral $\Sigma(\zeta;s)$ in (\ref{lem:frs asymptotic formula}), as $s\to\infty$. Observe that
$\Sigma(\zeta;s)=\Sigma_1(\zeta;s)+\Sigma_2(\zeta;s)$, where we have defined
\[\Sigma_i(\zeta;s):=\int_0^\epsilon e^{-s\cdot f(\zeta,(-1)^{i+1}\theta)}a(\zeta,(-1)^{i+1}\theta)d\theta\,, i=1,2.\]
Because of the similarity of $\Sigma_1(\zeta;s)$ and $\Sigma_2(\zeta;s)$, we analyze only the asymptotic behavior of $\Sigma_1(\zeta;s)$ under the hypotheses that $f(\zeta,\theta)$ and $a(\zeta,\theta)$ have respectively an $n$-to-$n$ and $p$-to-$q$ change of degree about $\theta=0$ as $\zeta\to \zeta_c$, and that $f(\zeta,\theta)$ has the properties stated in Lemma~\ref{lem:[zrws]F as an integral}. A similar analysis of the asymptotic behavior of $\Sigma_2(\zeta;s)$ is summarized at the end of this section.

Lemma~\ref{lem:[zrws]F as an integral} implies that $n\ge2$. In particular, we may write
 \begin{equation}\label{proof:expansion f}
 f(\zeta,\theta)=u(\zeta)\cdot\theta^n+\ldots
 \end{equation}
where $u$ is certain analytic function near $\zeta_c$ such that $u(\zeta_c)\ne0$. Since for all nonzero $\theta\in[-\epsilon,\epsilon]$, $\Re\{f(\zeta_c,\theta)\}>0$, we must have $\Re\{u(\zeta_c)\}\ge0$.
 
On the other hand, Lemma \ref{lem:p-to-q terms} implies that there is a unique representation of the form
\begin{equation}\label{proof:integral from 0 to theta of a(z,xi)}
\int_0^\theta a(\zeta,w)dw=\sum_{k=p}^q\frac{A_k(\zeta)}{k+1}\alpha^{k+1}\,,
\end{equation}
where $A_k(\zeta_c)=0$, for all $p\le k<q$, $A_q(\zeta_c)\ne0$, and $\alpha=\alpha(\zeta,\theta)$ is such that $\alpha(\zeta,0)=0$ and $\frac{\partial\alpha}{\partial\theta}(\zeta,0)=1$. The coefficients $A_k$, $p\le k\le q$, correspond to those appearing in Remark 1. The Inverse Mapping Theorem implies that $\Psi_1(\zeta,\theta):=(\zeta,\alpha(\zeta,\theta))$ is a biholomorphic map from an open neighborhood of $(\zeta,\theta)=(\zeta_c,0)$ to an open neighborhood of $(\zeta,\alpha)=(\zeta_c,0)$.  In particular, assuming that $\epsilon>0$ is sufficiently small, we can perform in $\Sigma_1(\zeta;s)$ the change of variables $\alpha=\alpha(\zeta,\theta)$ to obtain that
\[\Sigma_1(\zeta;s)=\sum_{k=p}^q A_k(\zeta)\int_0^{\alpha(\zeta,\epsilon)} e^{-s\cdot g(\zeta,\alpha)}\alpha^kd\alpha\,,\]
where $g(\zeta,\alpha):=f(\Psi_1^{-1}(\zeta,\alpha))$. This last function is analytic in a neighborhood of the origin. Furthermore, its Hartogs series about $(\zeta,\alpha)=(\zeta_c,0)$ in powers of $\alpha$ is of the form $g(\zeta,\alpha)=u(\zeta)\alpha^n+\ldots$ with $u(\zeta)$ is as in~(\ref{proof:expansion f}). This motivates us to consider the map
\[\Psi_2(\zeta,\alpha):=\left(\zeta,\alpha\cdot(u(\zeta))^{1/n}\cdot\left(1+\frac{g(\zeta,\alpha)-u(\zeta)\alpha^n}{u(\zeta)\alpha^n}\right)^{1/n}\right)\,,\]
where the principal branch of the $n$-th root function is to be used in both cases. Since $u(\zeta_c)\ne0$ and $\Re\{u(\zeta_c)\}\ge0$, it follows that $\Psi_2$ is well-defined and holomorphic near $(\zeta_c,0)$. Furthermore, if $\beta=\beta(\zeta,\alpha)$ is such that $\Psi_2(\zeta,\alpha)=(\zeta,\beta(\zeta,\alpha))$, the Inverse Mapping Theorem implies that $\Psi_2(\zeta,\alpha)$ is biholomorphic between open neighborhoods of  $(\zeta,\alpha)=(\zeta_c,0)$ and $(\zeta,\beta)=(\zeta_c,0)$. In particular, it follows that
$g(\Psi_2^{-1}(\zeta,\beta))=\beta^n$ and therefore
\[\Sigma_1(\zeta;s)=\sum_{k=p}^q A_k(\zeta)\int_0^{\beta(\zeta,\alpha(\zeta,\epsilon))} e^{-s\cdot \beta^n}(\alpha(\zeta,\beta))^k\frac{\partial\alpha}{\partial\beta}(\zeta,\beta)d\beta\,,\]
provided that $\epsilon>0$ is chosen sufficiently small to start with. We claim that the  domain of integration of the integrals participating in the summation above can be replaced by a real interval of the form $[0,\delta]$, for some $\delta>0$. For this observe that the condition $\Re\{f(\zeta_c,\epsilon)\}>0$ implies that $\Re\{(\beta(\zeta_c,\alpha(\zeta_c,\epsilon)))^n\}>0$. On the other hand, since $\beta(\zeta_c,\alpha(\zeta_c,\epsilon))=(u(\zeta_c))^{1/n}\epsilon+O(\epsilon^2)$, with $\Re\{u(\zeta_c)\}\ge0$, we conclude that $|\arg\{\beta(\zeta_c,\alpha(\zeta_c,\epsilon))\}|<\pi/(2n)$.
Since $\beta(\zeta,\alpha(\zeta,\epsilon))\to\beta(\zeta_c,\alpha(\zeta_c,\epsilon))$, as $\zeta\to \zeta_c$, we conclude that $|\arg\{\beta(\zeta,\alpha(\zeta,\epsilon))\}|<\pi/(2n)$, for all $\zeta$ sufficiently close to $\zeta_c$. Choosing $\delta:=\Re\{\beta(\zeta_c,\alpha(\zeta_c,\epsilon))\}$, it follows that there is a constant $c>0$ such that
\[\int_\delta^{\beta(\zeta,\alpha(\zeta,\epsilon))} e^{-s\cdot \beta^n}(\alpha(\zeta,\beta))^k\frac{\partial\alpha}{\partial\beta}(\zeta,\beta)d\beta=O(e^{-sc})\,,\]
as $s\to\infty$, uniformly for all $\zeta$ sufficiently close to $\zeta_c$ and for all $p\le k\le q$. This implies that
\begin{equation}\label{proof:Sigma1(z;s)}
\Sigma_1(\zeta;s)=\sum_{k=p}^q A_k(\zeta)\int_0^\delta e^{-s\cdot \beta^n}(\alpha(\zeta,\beta))^k\frac{\partial\alpha}{\partial\beta}(\zeta,\beta)d\beta+O(e^{-sc})\,,
\end{equation}
as $s\to\infty$, uniformly for all $\zeta$ sufficiently close to $\zeta_c$. An asymptotic expansion for the integrals participating in the summation above is easily obtained using the standard stationary phase method (see Chapter 6 in~\cite{BleHan86}). Indeed, since the Hartog's series of $(\alpha(\zeta,\beta))^k\frac{\partial\alpha}{\partial\beta}(\zeta,\beta)$ in powers of $\beta$ about $(\zeta,\beta)=(\zeta_c,0)$ must be of the form
\[(\alpha(\zeta,\beta))^k\frac{\partial\alpha}{\partial\beta}(\zeta,\beta)=\sum_{j=k}^\infty c_k(\zeta;j)\beta^j\,,\]
with $c_k(\zeta;k)=(u(\zeta))^{-(k+1)/n}$, then from (\ref{proof:Sigma1(z;s)}) it follows that
\begin{eqnarray}
\label{proof:asymptotic Sigma1(z;s)}
\Sigma_1(\zeta;s)&=&\sum_{k=p}^q A_k(\zeta)\cdot B_k(\zeta;s)+O(e^{-sc})\,,\\
\label{proof:asymptotic Bk(z;s)}
B_k(\zeta;s)&\approx&\sum_{j=k}^\infty\frac{c_k(\zeta;j)}{n}\Gamma\left(\frac{j+1}{n}\right)\cdot s^{-(j+1)/n}\,,
\end{eqnarray}
uniformly for all $\zeta$ sufficiently close to $\zeta_c$ as $s\to\infty$. The coefficients $c_k(\zeta;j)$ correspond to those appearing in Remark 2.  (\ref{proof:asymptotic Sigma1(z;s)}) and (\ref{proof:asymptotic Bk(z;s)}) provide a complete asymptotic description for $\Sigma_1(\zeta;s)$ which is uniform for all $\zeta$ sufficiently close to $\zeta_c$ as $s\to\infty$.

To obtain an asymptotic expansion for the term $\Sigma_2(\zeta;s)$, the uniqueness of the decomposition in (\ref{proof:integral from 0 to theta of a(z,xi)}) is relevant to relate the coefficients appearing in the expansion of $\Sigma_2(\zeta;s)$ with those in (\ref{proof:asymptotic Sigma1(z;s)}) and (\ref{proof:asymptotic Bk(z;s)}). Without delving into details it follows that
\begin{equation}
\label{proof:asymptotic Sigma2(z;s)}
\Sigma_2(\zeta;s)=\sum_{k=p}^q A_k(\zeta)\cdot \tilde{B}_k(\zeta;s)+O(e^{-sc})\,,
\end{equation}
where for the case in which $n$ is even it applies that
\begin{equation}
\tilde{B}_k(\zeta;s)\approx\sum_{j=k}^\infty\frac{(-1)^jc_k(\zeta;j)}{n}\Gamma\left(\frac{j+1}{n}\right)\cdot s^{-(j+1)/n}\,,
\end{equation}
however for the case in which $n$ is odd,
\begin{equation}\label{proof:n odd}
\tilde{B}_k(\zeta;s)\approx\sum_{j=k}^\infty\frac{(-1)^j D(j,n) c_k(\zeta;j)}{n}\Gamma\left(\frac{j+1}{n}\right)\cdot s^{-(j+1)/n}\,,
\end{equation}
where $D(j,n):=\exp\left(-\frac{i\pi(j+1)}{n}\cdot\hbox{sign}\{i[\theta^n]f(\zeta_c,\theta)\}\right)$. (\ref{thm:[zrws]F}) and (\ref{thm:Bk(z;s)}) in Theorem~\ref{thm:main} are now a direct consequence of (\ref{proof:asymptotic Sigma1(z;s)})-(\ref{proof:n odd}). This completes the proof of Theorem~\ref{thm:main}.$\Box$\\

\noindent{\bf Acknowledgments:} I would like to thank my graduate advisor, Robin Pemantle, and his collaborator, Mark Wilson, for their insights and support in my work on asymptotic analysis and generating functions. Special thanks also to Jean-Pierre Rosay and Saleh Tanveer for their helpful inputs in the preliminary version of my graduate dissertation which in one way or another have been reflected in here.

\end{document}